\PassOptionsToPackage{unicode}{hyperref}
\PassOptionsToPackage{hyphens}{url}
\documentclass[
  11pt,
]{article}
\usepackage{amsmath,amssymb}
\usepackage{iftex}
\ifPDFTeX
  \usepackage[T1]{fontenc}
  \usepackage[utf8]{inputenc}
  \usepackage{textcomp} 
\else 
  \usepackage{unicode-math} 
  \defaultfontfeatures{Scale=MatchLowercase}
  \defaultfontfeatures[\rmfamily]{Ligatures=TeX,Scale=1}
\fi
\usepackage{lmodern}
\usepackage{needspace}
\usepackage{titlesec}
\titleformat{\subsection}{\large\bfseries}{\thesubsection}{1em}{}
\titleformat{\subsubsection}{\normalsize\bfseries}{\thesubsubsection}{1em}{}
\ifPDFTeX\else
\fi
\IfFileExists{upquote.sty}{\usepackage{upquote}}{}
\IfFileExists{microtype.sty}{
  \usepackage[]{microtype}
  \UseMicrotypeSet[protrusion]{basicmath} 
}{}
\makeatletter
\@ifundefined{KOMAClassName}{
  \IfFileExists{parskip.sty}{%
    \usepackage{parskip}
  }{
    \setlength{\parindent}{0pt}
    \setlength{\parskip}{6pt plus 2pt minus 1pt}}
}{
  \KOMAoptions{parskip=half}}
\makeatother
\usepackage{xcolor}
\usepackage[margin=1in]{geometry}
\usepackage{longtable,booktabs,array}
\usepackage{calc} 
\usepackage{etoolbox}
\makeatletter
\patchcmd\longtable{\par}{\if@noskipsec\mbox{}\fi\par}{}{}
\makeatother
\IfFileExists{footnotehyper.sty}{\usepackage{footnotehyper}}{\usepackage{footnote}}
\makesavenoteenv{longtable}
\usepackage{graphicx}\usepackage{float}
\makeatletter
\def\maxwidth{\ifdim\Gin@nat@width>\linewidth\linewidth\else\Gin@nat@width\fi}
\def\maxheight{\ifdim\Gin@nat@height>\textheight\textheight\else\Gin@nat@height\fi}
\makeatother
\setkeys{Gin}{width=\maxwidth,height=\maxheight,keepaspectratio}
\makeatletter
\def\fps@figure{htbp}
\makeatother
\ifLuaTeX
\usepackage[bidi=basic]{babel}
\else
\usepackage[bidi=default]{babel}
\fi
\babelprovide[main,import]{american}

\def\languageshorthands#1{}
\ifLuaTeX
  \usepackage{selnolig}  
\fi
\IfFileExists{bookmark.sty}{\usepackage{bookmark}}{\usepackage{hyperref}}
\IfFileExists{xurl.sty}{\usepackage{xurl}}{} 
\hypersetup{
  pdftitle={Bounds for Orthogonal Expansion Coefficients via Weighted Primitives},
  pdfauthor={Vladimir Sluchak},
  pdflang={en-US},
  hidelinks,
  pdfcreator={LaTeX via pandoc}}

\title{Bounds for Orthogonal Expansion Coefficients via Weighted
Primitives}
\author{Vladimir Sluchak}
\date{}

\begin{document}
\maketitle

\hypertarget{abstract}{%
\subsection{Abstract}\label{abstract}}

For trigonometric and polynomial orthogonal expansions, we present
coefficient bounds with an explicit numerator determined by the total
variation of the function and the largest cumulative signed area of the
weighted basis element. This area, denoted by \(M_n\), is the sharp
constant in the corresponding weighted integral estimate. We calculate
or bound \(M_n\) for several classical systems, thereby obtaining
quantitative inequalities beyond an asymptotic statement of coefficient
decay. Integration between consecutive zeros gives a geometric
interpretation of the estimate. The comparison includes Chebyshev
polynomials of the first three kinds, Legendre, Jacobi, Gegenbauer, and
ordinary Laguerre polynomials, with normalization and boundary
conditions stated explicitly. Classical estimates and their elementary
consequences are organized as Fourier analysis in inequalities. In the
controlled-distortion application, ideal curves have selected
coefficients absent, and restrictions on the variation, corridor width,
and number of extrema of admissible departures give explicit guarantees
for those coefficients.

\hypertarget{introduction}{%
\subsection{1 Introduction}\label{introduction}}

An exact formula for a Fourier coefficient requires information about
the whole function. In applications, such information may be
unavailable, and the full curve may be impossible to prescribe. It is
then useful to replace exact coefficient calculation by inequalities
involving quantities that can be estimated or controlled, such as total
variation. In this sense, the problem is one of Fourier analysis in
inequalities.

The motivation comes from the control of selected harmonics of the
circular distribution of velocity upstream of a propeller. In the
formulation of {[}1{]}, ideal curves have selected harmonics absent, and
admissible departures from those curves are constrained so that the
corresponding harmonics remain small. The mathematical requirement is an
estimate connecting the variation of a curve, or of its departure from
an ideal curve, with its expansion coefficients. Section 10 develops
this use of the inequalities through explicit admissible corridors. The construction of ideal curves and their application to propeller inflow were introduced in [15] and subsequently discussed and illustrated in [16, 17]. A separate paper will present a fuller account and develop a broader formulation of the controlled-distortion method. We use the
term harmonics for the trigonometric components of the periodic inflow;
for a general orthogonal system we refer to expansion coefficients or
components.

The Riemann--Lebesgue theorem gives convergence of Fourier coefficients
to zero. Additional assumptions on variation or smoothness yield
quantitative decay estimates, including the classical \(1/n\) bound for
periodic functions of bounded variation; see Taibleson {[}4{]} and
Edwards {[}6{]}. Here bounded variation is the regularity assumption,
and the aim is to identify what features of a specified orthogonal
system permit analogous estimates. The explicit numerator is essential
for this purpose: an asymptotic order alone does not specify a design
tolerance at a selected harmonic.

Integration by parts gives the estimate \(|I_n(f)|\le V(f)M_n\). The
factor \(M_n\) is the smallest constant for which this inequality holds
for every function of bounded variation. It depends only on the chosen
basis function, its normalization, the weight, and the interval, and can
therefore be calculated independently of the function being expanded.
The zeros describe the oscillatory structure of the basis element, while
the weighted primitive measures the accumulated balance of its positive
and negative contributions.

The connection between bounded variation and primitives of general
orthonormal systems predates the present investigation. Tsagareishvili
{[}7{]} announced a necessary and sufficient condition for \(O(1/n)\)
coefficient decay for every BV function, under uniform boundedness and
zero-mean assumptions on the orthonormal system. His later paper {[}8,
Theorems 5 and 5'{]} gives grid and continuous primitive formulations.
The 1981 note states the results without proofs. On \([0,1]\), if
\(\|\phi_n\|_\infty\le B\), the grid quantity
\(K_n=\max_{0\le i\le n}|\int_0^{i/n}\phi_n(t)\,dt|\) satisfies

\[
K_n\le\left\|\int_0^x\phi_n(t)\,dt\right\|_\infty
\le K_n+\frac Bn.
\]

Thus the two primitive conditions have the same \(O(1/n)\) content under
those assumptions. These results concern bounds for each fixed BV
function, with the implied constant allowed to depend on that function.
Here the exact constant for a prescribed variation bound is kept
explicit. The zero-mean kernel argument also applies when the normalized
basis functions are not uniformly bounded. Gogoladze and Tsagareishvili
{[}5,9{]} provide further relevant results for general orthonormal
systems.

The integration by parts estimate is a classical consequence of the BV
bound for Stieltjes integrals; see, for example, Dragomir {[}10{]}.
Closely related Legendre coefficient and Gegenbauer primitive estimates
appear in Wang {[}11{]}, and weighted Gegenbauer estimates go back to
Durand {[}12{]}. Xiang and Liu {[}13, Theorem 1 with \(k=0\){]} give
precisely the Jacobi coefficient bound \(V(f)M_n/h_n\), with an explicit
weighted polynomial expression for \(M_n\). Their Remark 1 also locates
the primitive maxima at zeros of the original polynomial. Their uniform
weighted estimates yield the Jacobi orders compared here. Zhang, Xiang,
and Kong {[}14, Theorem 2.4{]} give coefficient decay for Laguerre
expansions of piecewise functions, including the \(O(n^{-3/4})\)
estimate for the ordinary Laguerre step-function example below. The
present treatment organizes these classical bounds and their
consequences into a unified comparison through \(M_n\), explains the
geometry by integration between zeros, and converts shape restrictions
into spectral guarantees in controlled distortion.

\hypertarget{notation-and-normalization}{%
\subsection{2 Notation and
normalization}\label{notation-and-normalization}}

Let \(I=[a,b]\) and let \(w\) be a nonnegative integrable weight.
Consider real orthogonal functions \(p_0,p_1,\ldots\) with \(p_0\)
constant, and write

\[
h_n=\|p_n\|_w^2=\int_a^b p_n(x)^2w(x)\,dx>0,
\qquad
I_n(f)=\int_a^b f(x)p_n(x)w(x)\,dx.
\]

The coordinate coefficient and the orthonormal coefficient are,
respectively,

\[
a_n(f)=\frac{I_n(f)}{h_n},
\qquad
b_n(f)=\frac{I_n(f)}{\sqrt{h_n}}.
\]

The bounds for individual coefficients do not require the system to
represent every function. Completeness is assumed only when we discuss
approximation by the full expansion; it means that finite linear
combinations of the basis functions are dense in \(L^2(w)\).

The total variation is

\[
V_I(f)=\sup_{a=x_0<\cdots<x_m=b}
\sum_{j=1}^{m}|f(x_j)-f(x_{j-1})|.
\]

We abbreviate this to \(V(f)\) when the interval is understood. If \(f\)
is absolutely continuous, then \(V(f)=\int_a^b|f'(x)|\,dx\). The
partition definition also covers jumps and singular continuous
variation. For a continuous curve with a finite partition into monotone
arcs, variation is the sum of the absolute changes along those arcs,
including the end arcs.

Throughout, \(n\ge1\). For orthogonal polynomials, \(p_n\) has degree
\(n\). Conventional normalizations are used unless stated otherwise.
Constants denoted by \(K\) may depend on the fixed parameters defining
the polynomial family, but are independent of the degree \(n\) and the
function \(f\). Their values may differ between estimates. For positive
sequences, \(u_n\asymp v_n\) means that there exist constants \(c,C>0\)
such that

\[
c\,v_n\le u_n\le C\,v_n
\]

for all sufficiently large \(n\). Thus the sequences have the same
asymptotic order.

\hypertarget{the-sharp-primitive-bound}{%
\subsection{3 The sharp primitive
bound}\label{the-sharp-primitive-bound}}

\textbf{Theorem 1 (classical BV estimate in primitive form).} Let \(k\)
be a real integrable function on \([a,b]\) such that \(\int_a^b k=0\),
and put

\[
G(x)=\int_a^x k(t)\,dt.
\]

For every real function \(f\) of bounded variation,

\[
\left|\int_a^b f(x)k(x)\,dx\right|
\le V(f)\|G\|_\infty.
\]

The constant is exact:

\[
\max_{\substack{f\in BV[a,b]\\V(f)\le1}}\left|\int_a^b f(x)k(x)\,dx\right|
=\|G\|_\infty.
\]

The maximum is attained in the BV class. The same value is the supremum
if \(f\) is required to be continuous and to have variation at most one;
in that smaller class it need not be attained. The threshold one is a
normalization: for any \(B>0\), the corresponding maximum over
\(V(f)\le B\) is \(B\|G\|_\infty\).

\emph{Proof.} The primitive \(G\) is absolutely continuous and satisfies
\(G(a)=G(b)=0\). Stieltjes integration by parts, valid here even when
\(f\) has jumps, gives

\[
\int_a^b f(x)k(x)\,dx
=-\int_a^b G(x)\,df(x).
\]

The absolute value of the last integral is at most
\(\|G\|_\infty V(f)\). To prove sharpness, suppose \(G\) is not
identically zero and choose \(c\in(a,b)\) with \(|G(c)|=\|G\|_\infty\).
The indicator function \(f=\mathbf{1}_{[c,b]}\), equal to zero before
\(c\) and one after \(c\), has variation one and satisfies

\[
\int_a^b fk=\int_c^b k=-G(c).
\]

Thus equality is attained in the BV class. Replacing the jump by a
monotone continuous ramp of width tending to zero preserves variation
one and makes the integrals converge to \(-G(c)\). The zero kernel case
is immediate. This proves the theorem.

Apply Theorem 1 to \(k=p_nw\). Orthogonality to constants gives zero
total integral. With

\[
G_n(x)=\int_a^x p_n(t)w(t)\,dt,
\qquad M_n=\|G_n\|_\infty,
\]

we obtain the sharp bounds

\[
|I_n(f)|\le V(f)M_n,
\qquad
|a_n(f)|\le\frac{V(f)M_n}{h_n},
\qquad
|b_n(f)|\le\frac{V(f)M_n}{\sqrt{h_n}}.
\]

In particular, a constant \(C\) satisfying

\[
|a_n(f)|\le\frac{C V(f)}{n h_n}
\]

for every \(n\ge1\) and every BV function exists if and only if
\(\sup_{n\ge1} nM_n<\infty\). For a fixed normalization, the best such
constant is \(\sup_n nM_n\). This statement concerns a bound uniform in
\(f\). It does not claim that every fixed function has the worst
possible coefficient behavior.

The same proof applies on \([0,\infty)\) when \(k\) is integrable with
zero integral, \(f\) has finite total variation there, and the primitive
satisfies \(G(0)=0\) and tends to zero at infinity. Such an \(f\) is
bounded. Integration by parts on finite intervals and passage to the
limit justify the formula.

\hypertarget{integration-between-zeros-and-cumulative-signed-areas}{%
\subsection{4 Integration between zeros and cumulative signed
areas}\label{integration-between-zeros-and-cumulative-signed-areas}}

For polynomials orthogonal with a positive weight on an interval, all
\(n\) zeros of the degree \(n\) polynomial are simple and lie inside the
interval {[}2, Section 18.2{]}. Thus the coefficient index is also the
number of sign changes. If

\[
a=z_0<z_1<\cdots<z_n<z_{n+1}=b,
\]

the zeros divide the interval into \(n+1\) sign intervals. The points
\(z_0\) and \(z_{n+1}\) are the endpoints, not additional zeros. Define,
for \(0\le k\le n\),

\[
A_k=\int_{z_k}^{z_{k+1}}p_n(x)w(x)\,dx,
\qquad
S_k=\sum_{r=0}^k A_r=G_n(z_{k+1}).
\]

On each interval \([z_k,z_{k+1}]\), the weighted polynomial \(p_nw\) has
a constant sign. For continuous \(f\), the weighted mean value theorem
therefore gives a point \(\xi_k\in[z_k,z_{k+1}]\) such that

\[
\int_{z_k}^{z_{k+1}} f(x)p_n(x)w(x)\,dx=f(\xi_k)A_k.
\]

These points are ordered, \(\xi_0\le\cdots\le\xi_n\). Summing the
interval integrals gives

\[
I_n(f)=\sum_{k=0}^n f(\xi_k)A_k.
\]

Since \(\sum A_k=0\), discrete summation by parts gives

\[
I_n(f)=\sum_{k=0}^{n-1}S_k
\bigl(f(\xi_k)-f(\xi_{k+1})\bigr).
\]

The variation bound enters through the ordered evaluation points:

\[
\sum_{k=0}^{n-1}|f(\xi_{k+1})-f(\xi_k)|\le V(f).
\]

Consequently,

\[
|I_n(f)|\le\max_{0\le k\le n-1}|S_k|
\sum_{k=0}^{n-1}|f(\xi_{k+1})-f(\xi_k)|
\le V(f)\max_{0\le k\le n-1}|S_k|.
\]

Between successive zeros, \(G_n\) is monotone. Since
\(G_n(a)=G_n(b)=0\), its largest absolute value is attained at an
interior zero of \(p_n\). Thus \(\max_{0\le k\le n-1}|S_k|=M_n\).
Division by \(h_n\) recovers the coefficient estimate of Theorem 1 for
continuous \(f\). The Stieltjes proof covers the whole BV class,
including discontinuous functions.

Equal absolute areas between zeros are unnecessary. The relevant
quantity is the largest cumulative signed area. Counting zeros describes
the oscillations, but does not alone bound this cumulative area.

\begin{figure}[H]
\centering
\includegraphics[width=0.95\textwidth,height=\textheight]{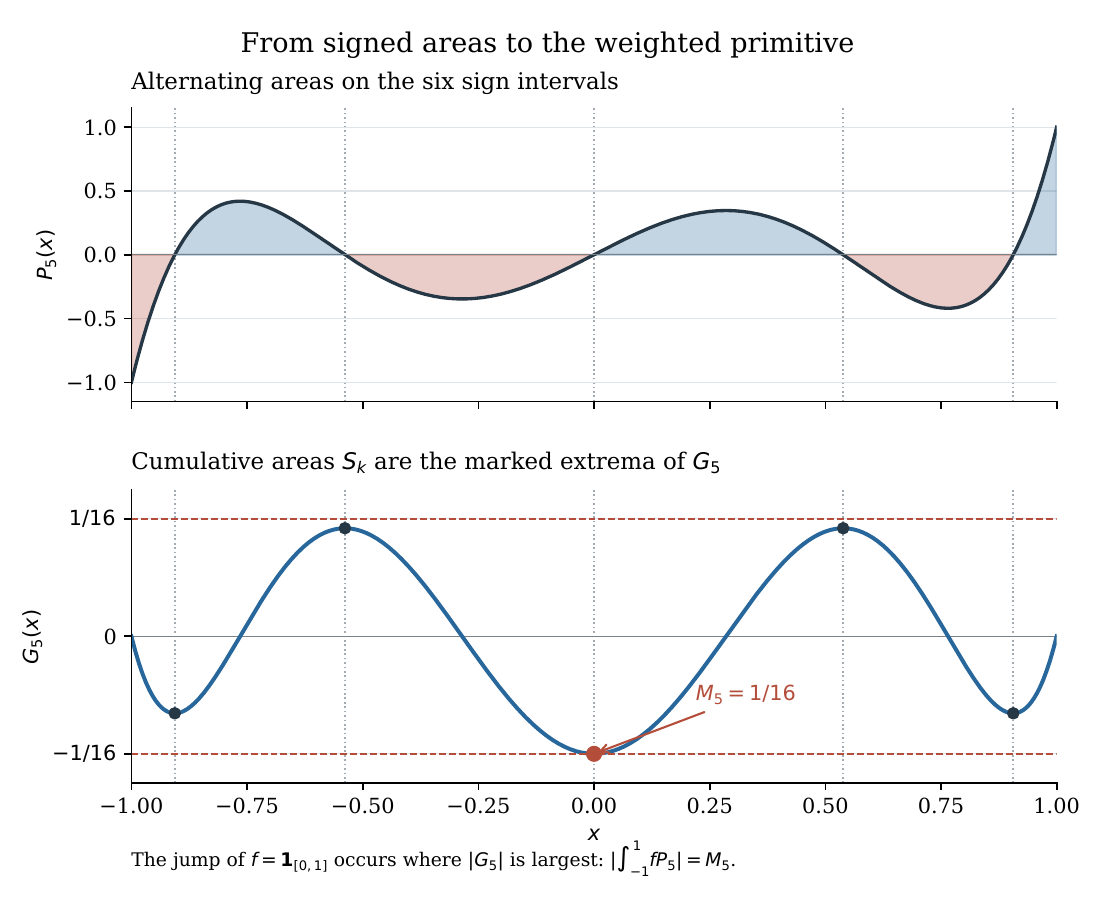}
\caption{Weighted primitive geometry for the Legendre polynomial of
degree five. The upper panel shows the alternating signed areas between
consecutive zeros. The lower panel shows their cumulative values
\(S_k=G_5(z_{k+1})\). The maximum absolute value is \(M_5=1/16\),
attained at the central zero. A step with its jump there attains the
integral bound.}
\end{figure}

\hypertarget{elementary-family-estimates}{%
\subsection{5 Elementary family
estimates}\label{elementary-family-estimates}}

\hypertarget{the-cosine-system-and-periodic-trigonometric-functions}{%
\subsubsection{5.1 The cosine system and periodic trigonometric
functions}\label{the-cosine-system-and-periodic-trigonometric-functions}}

For \(p_n(x)=\cos(nx)\) on \([0,\pi]\) with \(w=1\),

\[
G_n(x)=\frac{\sin(nx)}n,
\qquad M_n=\frac1n,
\qquad h_n=\frac\pi2.
\]

Consequently \(|a_n(f)|\le2V(f)/(\pi n)\), with the exact integral
constant \(1/n\).

On a full period \([0,2\pi]\), assume \(f(0)=f(2\pi)\). The primitives
\(\sin(nx)/n\) and \(-\cos(nx)/n\) give, by integration by parts,

\[
\left|\int_0^{2\pi}f(x)\cos(nx)\,dx\right|\le\frac{V(f)}n,
\qquad
\left|\int_0^{2\pi}f(x)\sin(nx)\,dx\right|\le\frac{V(f)}n.
\]

The sine boundary term vanishes because the endpoint values of \(f\)
agree. Both coordinate coefficients, whose squared basis norms are
\(\pi\), are bounded by \(V(f)/(\pi n)\). This is the periodic
convention used in the motivating estimate {[}1{]}. Without periodic
endpoint matching, the anchored sine primitive \((1-\cos(nx))/n\)
instead gives a bound \(2V(f)/n\) for the unnormalized integral. The
boundary convention therefore belongs in the statement.

\hypertarget{chebyshev-polynomials-of-the-first-kind}{%
\subsubsection{5.2 Chebyshev polynomials of the first
kind}\label{chebyshev-polynomials-of-the-first-kind}}

Let \(p_n=T_n\) and \(w(x)=(1-x^2)^{-1/2}\) on \([-1,1]\). With
\(x=\cos\theta\),

\[
T_n(\cos\theta)=\cos(n\theta),
\qquad
G_n(\cos\theta)=-\frac{\sin(n\theta)}n.
\]

Hence \(M_n=1/n\). Since \(h_n=\pi/2\),

\[
|a_n(f)|\le\frac{2V(f)}{\pi n}.
\]

The substitution is monotone, so
\(V_{[0,\pi]}(f(\cos\theta))=V_{[-1,1]}(f)\). The result is the cosine
estimate in a different variable.

\hypertarget{legendre-polynomials}{%
\subsubsection{5.3 Legendre polynomials}\label{legendre-polynomials}}

For \(p_n=P_n\) and \(w=1\), the standard identities in {[}2, Section
18.9{]} give

\[
G_n(x)=\frac{P_{n+1}(x)-P_{n-1}(x)}{2n+1}.
\]

At any zero \(z\) of \(P_n\), the recurrence relation implies
\((n+1)P_{n+1}(z)=-nP_{n-1}(z)\). Therefore

\[
G_n(z)=-\frac{P_{n-1}(z)}{n+1}.
\]

Using \(|P_{n-1}(x)|\le1\) on \([-1,1]\), and the location of the
extrema of \(G_n\), we obtain

\[
M_n\le\frac1{n+1},
\qquad
|a_n(f)|\le\frac{V(f)}{(n+1)h_n},
\qquad h_n=\frac2{2n+1}.
\]

This convenient explicit bound is stronger than \(V(f)/(nh_n)\), but is
not asymptotically optimal. Section 6 gives \(M_n=O(n^{-3/2})\). Wang
{[}11, Theorem 3.1, case \(m=0\){]} uses the same Stieltjes argument and
a sharper bound for \(P_{n+1}-P_{n-1}\) to obtain an explicit
\(O(V(f)n^{-1/2})\) coordinate coefficient estimate. Our elementary
estimate is included for its short zero-based proof, rather than as an
improvement on that result.

\hypertarget{chebyshev-polynomials-of-the-second-kind}{%
\subsubsection{5.4 Chebyshev polynomials of the second
kind}\label{chebyshev-polynomials-of-the-second-kind}}

For \(p_n=U_n\) and \(w(x)=\sqrt{1-x^2}\), use

\[
U_n(\cos\theta)=\frac{\sin((n+1)\theta)}{\sin\theta}.
\]

Product to sum gives

\[
G_n(\cos\theta)=\int_\theta^\pi\sin((n+1)u)\sin u\,du
=-\frac{\sin(n\theta)}{2n}
+\frac{\sin((n+2)\theta)}{2(n+2)}.
\]

Consequently

\[
M_n\le\frac1{2n}+\frac1{2(n+2)}
=\frac{n+1}{n(n+2)}<\frac1n.
\]

With \(h_n=\pi/2\), this gives \(|a_n(f)|\le V(f)(n+1)/(n(n+2)h_n)\).

\hypertarget{chebyshev-polynomials-of-the-third-kind}{%
\subsubsection{5.5 Chebyshev polynomials of the third
kind}\label{chebyshev-polynomials-of-the-third-kind}}

Use the convention

\[
V_n(\cos\theta)=\frac{\cos((n+\tfrac12)\theta)}{\cos(\theta/2)},
\qquad
w(x)=\sqrt{\frac{1+x}{1-x}}.
\]

Here \(h_n=\pi\). The transformed integrand is
\(2\cos(u/2)\cos((n+\tfrac12)u)=\cos(nu)+\cos((n+1)u)\), and hence

\[
G_n(\cos\theta)
=-\frac{\sin(n\theta)}n-\frac{\sin((n+1)\theta)}{n+1}.
\]

Thus \(M_n\le1/n+1/(n+1)<2/n\), giving \(|a_n(f)|\le2V(f)/(nh_n)\). The
best constant uniform in \(n\) is exactly \(2\). Indeed, at
\(\theta_n=\pi/(2n+1)\),

\[
n|G_n(\cos\theta_n)|
=\left(1+\frac n{n+1}\right)
\cos\left(\frac\pi{4n+2}\right)\longrightarrow2.
\]

Together with the upper bound, this proves \(nM_n\to2\). The fourth
kind, with the reflected convention and reflected weight, has the same
primitive norms by the change \(x\mapsto-x\).

These bounds can also be recovered directly from classical cosine
coefficient estimates. Set \(F(\theta)=f(\cos\theta)\) and
\(c_j=(2/\pi)\int_0^\pi F(\theta)\cos(j\theta)\,d\theta\). Product
identities give, for \(n\ge1\),

\[
a_n^{\mathrm{II}}=\frac{c_n-c_{n+2}}2,
\qquad
a_n^{\mathrm{III}}=\frac{c_n+c_{n+1}}2.
\]

Applying \(|c_j|\le2V(f)/(\pi j)\) gives the same upper bounds. The
weighted primitives make the associated sharp integral constant \(M_n\)
visible; the displayed upper inequalities themselves are elementary
consequences of the classical cosine estimate.

\hypertarget{jacobi-polynomials}{%
\subsection{6 Jacobi polynomials}\label{jacobi-polynomials}}

Let \(p_n=P_n^{(\alpha,\beta)}\) have the conventional Jacobi
normalization and weight

\[
w(x)=(1-x)^\alpha(1+x)^\beta,
\qquad \alpha,\beta>-1.
\]

The weighted derivative identity {[}2, formula 18.9.16{]}, with its
parameters shifted by one, gives

\[
\frac{d}{dx}\left[(1-x)^{\alpha+1}(1+x)^{\beta+1}
P_{n-1}^{(\alpha+1,\beta+1)}(x)\right]
=-2n w(x)P_n^{(\alpha,\beta)}(x).
\]

The bracketed function vanishes at both endpoints, so

\[
G_n(x)=-\frac1{2n}(1-x)^{\alpha+1}(1+x)^{\beta+1}
P_{n-1}^{(\alpha+1,\beta+1)}(x).
\]

\textbf{Proposition 2.} For fixed \(\alpha,\beta>-1\),

\[
M_n\le K_{\alpha,\beta}n^{-q},
\qquad q=\min\left\{\frac32,\alpha+2,\beta+2\right\}>1.
\]

In particular, \(M_n=O(1/n)\) throughout this parameter range.

If \(\alpha,\beta\ge-1/2\), then \(M_n=O(n^{-3/2})\).

\emph{Proof.} We use the standard uniform Jacobi envelope, for fixed
\(A,B>0\),

\[
|P_{n-1}^{(A,B)}(\cos\theta)|
\le K n^{-1/2}
(\theta+n^{-1})^{-A-1/2}
(\pi-\theta+n^{-1})^{-B-1/2}.
\]

This envelope follows from the uniform endpoint Bessel expansions in
{[}2, formulas 18.15.6--18.15.8{]}, together with their reflected
versions and the interior Darboux expansion. It is a classical input to
the proof, rather than an estimate derived here from orthogonality
alone. To see the endpoint scaling, the small argument bound
\(J_A(u)=O(u^A)\) gives \(O(n^A)\) for the polynomial when
\(\theta\le1/n\); for \(\theta\ge1/n\), the bound \(J_A(u)=O(u^{-1/2})\)
gives the factor \(n^{-1/2}\theta^{-A-1/2}\) away from the opposite
endpoint. Fixed small degrees are absorbed into \(K\).

Set \(A=\alpha+1\), \(B=\beta+1\). On \(0\le\theta\le\pi/2\), the
opposite endpoint factors are bounded. Since \(1-\cos\theta\) is
comparable to \(\theta^2\),

\[
|G_n(\cos\theta)|\le
K n^{-3/2}\theta^{2\alpha+2}
(\theta+n^{-1})^{-\alpha-3/2}.
\]

For \(0\le\theta\le1/n\), the right hand side is at most
\(K n^{-\alpha-2}\). For \(1/n\le\theta\le\pi/2\), it is bounded by
\(K n^{-3/2}\theta^{\alpha+1/2}\). Its maximum there is \(O(n^{-3/2})\)
when \(\alpha\ge-1/2\), and \(O(n^{-\alpha-2})\) otherwise. Reflection
gives the corresponding bounds involving \(\beta\) on the other half of
the interval. Combining the bounds proves the proposition.

The coefficient estimate is therefore

\[
|a_n(f)|\le\frac{K_{\alpha,\beta}V(f)}{n^q h_n}.
\]

The conventional Jacobi squared norm satisfies \(h_n\asymp n^{-1}\) for
fixed parameters. Accordingly, this is a bound \(O(V(f)n^{1-q})\) for
coordinate coefficients and \(O(V(f)n^{1/2-q})\) for orthonormal
coefficients. For \(\alpha,\beta\ge-1/2\), these become
\(O(V(f)n^{-1/2})\) and \(O(V(f)n^{-1})\), respectively. The primitive
rate must not be confused with the coordinate coefficient rate.

These estimates are consistent with the explicit BV bounds of Xiang and
Liu {[}13{]}. Their Theorem 1, Eq. (14), at \(k=0\) is exactly
\(|a_n(f)|\le V(f)M_n/h_n\) with the primitive displayed above. Their
Theorem 4, Eq. (34), applied with parameters \(\alpha+1,\beta+1>0\),
bounds the weighted polynomial by
\(O(n^{-\min\{1/2,\alpha+1,\beta+1\}})\) and hence gives the same
exponent \(q\) for \(M_n\). Their Remark 1, Eq. (17), identifies the
zeros at which the maximum can be computed. The proof here spells out
the endpoint and interior contributions in the notation used throughout
this comparison.

\hypertarget{gegenbauer-normalization}{%
\subsection{7 Gegenbauer normalization}\label{gegenbauer-normalization}}

Restrict here to \(\lambda>0\). Conventional Gegenbauer polynomials
satisfy {[}2, formula 18.7.1{]}

\[
C_n^{(\lambda)}(x)=d_n P_n^{(\lambda-1/2,\lambda-1/2)}(x),
\qquad
d_n=\frac{(2\lambda)_n}{(\lambda+\tfrac12)_n}
\asymp n^{\lambda-1/2},
\]

where \((c)_n\) denotes the rising factorial. Both systems use the
weight \((1-x^2)^{\lambda-1/2}\), so their primitives and their uniform
norms satisfy

\[
G_n^{\mathrm{Geg}}=d_nG_n^{\mathrm{Jac}},
\qquad M_n^{\mathrm{Geg}}=d_nM_n^{\mathrm{Jac}}.
\]

Since the corresponding Jacobi parameters exceed \(-1/2\), Proposition 2
gives \(M_n^{\mathrm{Jac}}=O(n^{-3/2})\). Multiplication by
\(d_n\asymp n^{\lambda-1/2}\) therefore gives

\[
M_n^{\mathrm{Geg}}=O(n^{\lambda-2}).
\]

This order also has a matching lower bound. On a closed interior
interval, the Darboux formula for the shifted Jacobi polynomial in the
primitive has the form
\(n^{-1/2}A(\theta)\cos(n\theta+\gamma(\theta))+O(n^{-3/2})\), where
\(A\) is positive and bounded away from zero. The phase is affine in
\(\theta\) with slope \(n+O(1)\). For every sufficiently large \(n\), a
point in that interval can be chosen where the cosine has magnitude one.
The weight factors remain bounded away from zero. Multiplication by
\(d_n/n\) therefore gives a lower bound of order \(n^{\lambda-2}\). Thus

\[
M_n^{\mathrm{Geg}}\asymp n^{\lambda-2}.
\]

It follows that the conventional primitive condition \(M_n=O(1/n)\)
holds exactly for \(0<\lambda\le1\). Legendre and Chebyshev II occur at
\(\lambda=1/2\) and \(\lambda=1\).

The same primitive can be written directly as

\[
G_n^{\mathrm{Geg}}(x)=
-\frac{2\lambda}{n(n+2\lambda)}(1-x^2)^{\lambda+1/2}
C_{n-1}^{(\lambda+1)}(x).
\]

Durand's weighted Gegenbauer bound {[}12, Eq. (19){]} gives the explicit
consequence

\[
M_n^{\mathrm{Geg}}\le
\frac{2\Gamma((n+1)/2+\lambda)}
{n(n+2\lambda)\Gamma(\lambda)\Gamma((n+1)/2)}.
\]

The ratio of gamma functions is of order \(n^\lambda\), which again
gives \(O(n^{\lambda-2})\). Wang {[}11, Section 4.1{]} studies the
corresponding generalized Gegenbauer--Gauss--Lobatto function, equal to
\(2(n+\lambda)G_n^{\mathrm{Geg}}\), and gives detailed maximum
estimates. These results are direct predecessors of the primitive
comparison used here.

The same Durand estimate is Lemma 2, Eq. (21), of Xiang and Liu
{[}13{]}. Their Section 2.1 uses it to give explicit Gegenbauer
coefficient bounds and extensions to derivatives of bounded variation.
For the BV case considered here, their Jacobi Theorem 1 together with
this lemma gives the displayed primitive bound after changing
normalization. Their Theorem 3 is stated for derivative orders
\(k\ge1\); no extension of that theorem's stated assumptions is needed
for our \(k=0\) comparison.

An explicit example makes the failure for larger parameters visible. At
\(\lambda=2\),

\[
G_n(x)=-\frac4{n(n+4)}(1-x^2)^{5/2}C_{n-1}^{(3)}(x).
\]

For \(n=2m+1\), using \(C_{2m}^{(3)}(0)=(-1)^m(m+1)(m+2)/2\),

\[
|G_{2m+1}(0)|=
\frac{2(m+1)(m+2)}{(2m+1)(2m+5)}\longrightarrow\frac12.
\]

These conventional primitives do not tend uniformly to zero. Compactness
of the orthogonality interval alone is therefore insufficient for the
primitive condition.

For orthonormal Gegenbauer polynomials the result is different. Since
\(h_n\asymp n^{2\lambda-2}\), their primitives have norm

\[
\frac{M_n^{\mathrm{Geg}}}{\sqrt{h_n}}\asymp\frac1n
\qquad(\lambda>0).
\]

This agrees with the corresponding orthonormal Jacobi system. Rescaling
changes \(M_n\) and the coordinate coefficients, while the ratio
\(M_n/\sqrt{h_n}\) is invariant under a nonzero rescaling of each basis
member. Ordering still matters and is fixed here by degree. Parameters
\(-1/2<\lambda<0\) are outside the scope of this section.

\hypertarget{ordinary-laguerre-polynomials}{%
\subsection{8 Ordinary Laguerre
polynomials}\label{ordinary-laguerre-polynomials}}

Let \(p_n=L_n\) on \([0,\infty)\) with weight \(e^{-x}\). Ordinary
Laguerre polynomials are already orthonormal for this weight: \(h_n=1\).
Formula 18.9.24 in {[}2{]} yields

\[
G_n(x)=\frac1n e^{-x}xL_{n-1}^{(1)}(x).
\]

For each degree, this vanishes at zero and tends to zero at infinity. We
show that its uniform norm is not \(O(1/n)\).

For fixed \(x>0\), the classical Laguerre expansion {[}2, formulas
18.15.14--18.15.16{]} gives

\[
G_n(x)=\frac{e^{-x/2}x^{1/4}}{\sqrt\pi}n^{-3/4}
\left[
\cos\left(2\sqrt{(n-1)x}-\frac{3\pi}4\right)
+O(n^{-1/2})\right].
\]

The phase increases without bound and its increments tend to zero. For
each sufficiently large multiple of \(2\pi\), select the first index for
which the phase exceeds it. The overshoot tends to zero. Along the
resulting subsequence the cosine tends to one, so
\(|G_n(x)|\ge c_x n^{-3/4}\). Consequently \(nM_n\) is unbounded.

This gives a fixed BV counterexample to \(1/n\) coefficient decay. Take
\(f(t)=\mathbf{1}_{[x,\infty)}(t)\), with \(V(f)=1\). Orthogonality to
constants and \(h_n=1\) imply \(a_n(f)=-G_n(x)\), which is not
\(O(1/n)\). Continuous ramps also rule out a \(1/n\) inequality uniform
over continuous functions of variation at most one, by Theorem 1. The
fixed step example alone does not specify the optimal rate for each
individual continuous function.

This step belongs to the class \(z(t)(t-a)_+^\mu\) studied by Zhang,
Xiang, and Kong {[}14{]}: take \(z\equiv1\), \(\mu=0\), and \(a=x>0\).
Their convention assigns the midpoint value \(1/2\) at the jump;
changing a value at one point leaves the coefficients unchanged. Their
norm factor \(\Gamma(n+\alpha+1)/n!\) also reduces to one when
\(\alpha=0\). Their Theorem 2.4, p.~6, gives, for the class defined in
their Eq. (1),

\[
|a_n^{(\alpha)}(f)|=O\left(n^{-(\alpha+\mu)/2-3/4}\right).
\]

In our step case this is \(O(n^{-3/4})\), in agreement with the direct
primitive calculation. The oscillatory formula above supplies the
subsequence on which the \(1/n\) bound fails. Thus this example
illustrates an established decay phenomenon through a short exact
primitive identity. The rate at a fixed jump location does not by itself
determine \(M_n\), which takes the supremum over all locations for each
degree.

This argument establishes failure for ordinary Laguerre polynomials. It
neither gives a global upper estimate \(M_n=O(n^{-3/4})\) nor classifies
all generalized Laguerre parameters; both would require additional
analysis.

\hypertarget{comparison-of-the-primitive-estimates}{%
\subsection{9 Comparison of the primitive
estimates}\label{comparison-of-the-primitive-estimates}}

The bounds below concern \(M_n\), the sharp constant for the
unnormalized weighted integral. To obtain a coordinate coefficient
bound, multiply by \(V(f)/h_n\).

\begin{longtable}[]{@{}
  >{\raggedright\arraybackslash}p{(\columnwidth - 2\tabcolsep) * \real{0.5000}}
  >{\raggedright\arraybackslash}p{(\columnwidth - 2\tabcolsep) * \real{0.5000}}@{}}
\toprule\noalign{}
\begin{minipage}[b]{\linewidth}\raggedright
Conventional system
\end{minipage} & \begin{minipage}[b]{\linewidth}\raggedright
Estimate for \(M_n\)
\end{minipage} \\
\midrule\noalign{}
\endhead
\bottomrule\noalign{}
\endlastfoot
Cosine on \([0,\pi]\) & \(M_n=1/n\) \\
Chebyshev I & \(M_n=1/n\) \\
Legendre & \(M_n\le1/(n+1)\) and \(M_n=O(n^{-3/2})\) \\
Chebyshev II & \(M_n\le(n+1)/(n(n+2))<1/n\) \\
Chebyshev III & \(M_n<2/n\) and \(nM_n\to2\) \\
Jacobi with \(\alpha,\beta>-1\) & \(M_n=O(n^{-q})\),
\(q=\min\{3/2,\alpha+2,\beta+2\}>1\) \\
Gegenbauer with \(\lambda>0\) & \(M_n\asymp n^{\lambda-2}\) \\
Ordinary Laguerre & \(M_n\) is not \(O(1/n)\) \\
\end{longtable}

The general coefficient inequality is sharp with \(M_n\) itself. A
convenient upper estimate for \(M_n\), such as \(1/(n+1)\) in the
Legendre case, need not be the best numerical constant for each degree.

The numerical illustration below distinguishes a primitive norm from a
fixed-point witness. In its upper panel the extrema of each primitive
are evaluated at the polynomial zeros. In its lower panel only
\(|G_n(1)|\) is plotted for Laguerre, not the global maximum \(M_n\).
Since \(|G_n(1)|\le M_n\), the growing subsequence of \(n|G_n(1)|\)
already rules out boundedness of \(nM_n\). The calculations illustrate
the analytic arguments; they are not used as proofs.

\begin{figure}[H]
\centering
\includegraphics[width=0.95\textwidth,height=\textheight]{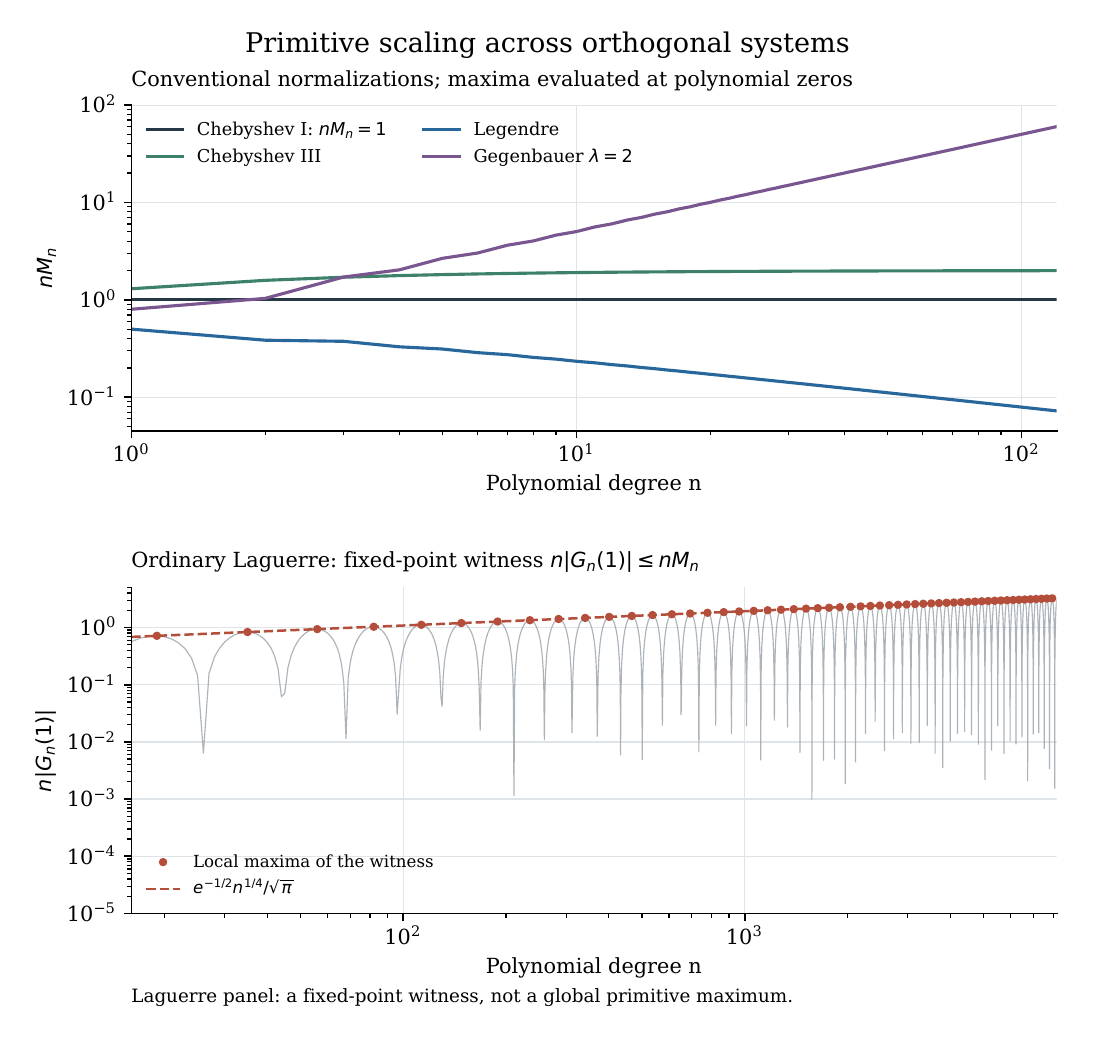}
\caption{Primitive scaling in conventional normalizations. Upper panel:
\(nM_n\) for Chebyshev I, Chebyshev III, Legendre, and Gegenbauer with
\(\lambda=2\); each maximum is evaluated at the zeros of its basis
polynomial. Lower panel: the Laguerre fixed-point witness \(n|G_n(1)|\),
its local maxima, and the asymptotic envelope
\(e^{-1/2}n^{1/4}/\sqrt\pi\). This lower panel does not assert a global
estimate for \(M_n\).}
\end{figure}

\hypertarget{controlled-distortion-and-other-consequences}{%
\subsection{10 Controlled distortion and other
consequences}\label{controlled-distortion-and-other-consequences}}

\hypertarget{bounds-from-incomplete-information}{%
\subsubsection{10.1 Bounds from incomplete
information}\label{bounds-from-incomplete-information}}

Suppose only \(V(f)\le B\) is known. Theorem 1 gives
\(|a_n(f)|\le B M_n/h_n\) for every admissible function, without
calculating its coefficients. Sharpness means that no smaller bound can
be guaranteed from this variation information alone. A smaller bound may
become possible if additional restrictions, such as monotonicity on
prescribed subintervals or boundary conditions, are imposed.

\hypertarget{ideal-curves-and-variation-corridors}{%
\subsubsection{10.2 Ideal curves and variation
corridors}\label{ideal-curves-and-variation-corridors}}

Let \(f_*\) be an ideal curve and set \(u=f-f_*\). Linearity gives

\[
|a_n(f)-a_n(f_*)|\le\frac{V(u)M_n}{h_n}.
\]

Fix a finite nonempty target set \(\mathcal T\) and suppose
\(a_n(f_*)=0\) for every \(n\in\mathcal T\). Given positive tolerances
\(\varepsilon_n\), put

\[
\eta=\min_{n\in\mathcal T}\frac{\varepsilon_n h_n}{M_n}.
\]

Every curve in the admissible class

\[
\mathcal C_\eta(f_*)=\{f:V(f-f_*)\le\eta\}
\]

satisfies \(|a_n(f)|\le\varepsilon_n\) simultaneously for all
\(n\in\mathcal T\). The ideal itself need not be physically attainable:
the purpose is to identify attainable members of a class with guaranteed
spectral behavior. Since constants do not affect the target coefficients
or total variation, a prescribed mean may be imposed separately. If the
target coefficients of \(f_*\) are small rather than zero, the
corresponding guarantee is \(|a_n(f)|\le|a_n(f_*)|+\eta M_n/h_n\).

This is the mathematical use of the inequalities in controlled
distortion {[}1{]}. Total variation is applied to the departure from the
ideal, rather than only to the original flow curve. The result supplies
sufficient spectral guarantees; whether a curve in this class can be
produced by a realizable flow modification is a separate design
question.

\hypertarget{geometric-corridors-and-restrictions-on-extrema}{%
\subsubsection{10.3 Geometric corridors and restrictions on
extrema}\label{geometric-corridors-and-restrictions-on-extrema}}

A bound on maximum pointwise departure alone does not bound total
variation. One useful additional restriction, employed in the
controlled-distortion construction of {[}1{]}, is a bound on the number
of extrema of the departure \(u=f-f_*\). Suppose \(u\) is continuous and
piecewise monotone, \(|u(x)|\le\delta\), and there are at most \(r\)
interior turning points, with monotonicity between successive turning
points and the endpoints. There are then at most \(r+1\) monotone arcs.
Each contributes at most \(2\delta\) to the variation, so

\[
V(u)\le2\delta(r+1),
\qquad
|a_n(f)-a_n(f_*)|\le
\frac{2\delta(r+1)M_n}{h_n}.
\]

For a nonconstant periodic departure with \(r\) alternating extrema
counted around the period and monotonicity between successive extrema,
there are \(r\) cyclic monotone arcs. In that case,

\[
V(u)\le2\delta r=r d,
\qquad d=2\delta,
\]

where \(d\) is the full width of the corridor around the ideal curve. A
constant departure has variation zero. These relations connect corridor
width and the number of extrema with a guaranteed bound on selected
coefficients. The extrema to be counted are those of \(f-f_*\), not
necessarily those of \(f\).

For periodic inflow velocity \(V(\varphi)\), let \(V_*(\varphi)\) have
zero sine and cosine coefficients at the target indices. If the periodic
departure \(u=V-V_*\) lies in a corridor of full width \(d\) and has at
most \(r\) cyclic extrema as above, Section 5.1 gives, for the real
Fourier coefficients,

\[
|A_n(V)|\le\frac{rd}{\pi n},
\qquad
|B_n(V)|\le\frac{rd}{\pi n}
\qquad(n\in\mathcal T).
\]

Here \(A_n\) and \(B_n\) denote the cosine and sine coordinate
coefficients in the usual real Fourier convention, with factor
\(1/\pi\). Write the amplitude of the \(n\)th real harmonic as
\(H_n=(A_n^2+B_n^2)^{1/2}\). In fact the same bound holds for this
amplitude:

\[
H_n(V)\le\frac{V(u)}{\pi n}\le\frac{rd}{\pi n}.
\]

Indeed, for every phase \(\psi\), the coefficient
\(A_n\cos\psi+B_n\sin\psi\) is the integral against
\(\cos(n\varphi-\psi)\) divided by \(\pi\). Its periodic primitive has
absolute value at most \(1/n\), and maximizing over \(\psi\) gives
\(H_n\). Thus \(d\le\min_{n\in\mathcal T}(\pi n\varepsilon_n/r)\) is a
sufficient corridor-width condition for \(H_n(V)\le\varepsilon_n\) when
\(r>0\).

\hypertarget{choices-of-an-ideal-curve}{%
\subsubsection{10.4 Choices of an ideal
curve}\label{choices-of-an-ideal-curve}}

The reference curve can be chosen to suit the design problem. A constant
ideal removes every nonzero harmonic. A dynamic ideal retains a selected
small set of major harmonics of the original curve, provided the target
harmonics are excluded. An optimal ideal removes only the target sine
and cosine components; it is the orthogonal projection onto the subspace
where those components vanish and minimizes the unweighted mean-square
departure on the period. A partial ideal coincides with the original
outside prescribed sectors and satisfies the target harmonic conditions
through modifications within those sectors. Modifications of peaks form
a subclass of the partial construction. Existence under additional
amplitude, sign, or physical restrictions must be examined in each
application.

For all these choices, the same coefficient inequality turns
restrictions on the departure into sufficient spectral guarantees. The construction of these ideal curves and their application to the propeller-flow problem were introduced in [15] and subsequently discussed and illustrated in [16, 17]. A separate paper will present a fuller account and develop a broader formulation of the controlled-distortion method. The
present application isolates the mathematical step shared by them.

\hypertarget{truncated-expansions}{%
\subsubsection{10.5 Truncated expansions}\label{truncated-expansions}}

Assume the orthogonal system is complete in \(L^2(w)\) and
\(f\in L^2(w)\). For \(S_Nf=\sum_{n=0}^N a_n(f)p_n\), Parseval's
identity gives

\[
\|f-S_Nf\|_w^2
=\sum_{n>N}|a_n(f)|^2h_n
\le V(f)^2\sum_{n>N}\frac{M_n^2}{h_n}.
\]

The bound is useful when the right hand series converges. For
conventional Jacobi polynomials, Proposition 2 and \(h_n\asymp n^{-1}\)
give

\[
\|f-S_Nf\|_w\le K_{\alpha,\beta}V(f)N^{1-q},
\qquad N\ge1.
\]

Indeed, \(\sum_{n>N}n^{1-2q}=O(N^{2-2q})\) because \(q>1\). In
particular, the bound is \(O(V(f)N^{-1/2})\) when
\(\alpha,\beta\ge-1/2\). These estimates provide a quantitative
approximation guarantee, rather than a claim of an optimal error rate
for every fixed function.

\hypertarget{scope-and-further-questions}{%
\subsection{11 Scope and further
questions}\label{scope-and-further-questions}}

The general principle applies to any real integrable zero mean kernel,
so orthogonality is sufficient for producing the kernels but is not
required by Theorem 1. Orthogonality becomes important in identifying
coefficients, interpreting their indices, and relating their squares to
approximation error.

The remaining classification problem must specify the functions, their
ordering, normalization, weight, and boundary conditions. One may study
either \(M_n\) in a prescribed conventional normalization or the
normalization invariant quantity \(M_n/\sqrt{h_n}\) for orthonormal
coefficients. The latter permits direct comparisons between differently
normalized descriptions of the same orthogonal system.

Further work may examine optimal constants for individual degrees,
generalized Laguerre parameters, and systems of eigenfunctions arising
from physical boundary value problems. For such systems, the natural
index may be related to an eigenvalue or spatial scale rather than
polynomial degree. The question remains concrete: how does the largest
cumulative signed area depend on that natural index, and what
coefficient guarantees follow for functions with bounded variation?

The family estimates presented here are established results or
consequences of the classical identities and bounds cited above. In
particular, {[}13{]} supplies the Jacobi BV coefficient inequality and
the weighted estimates behind its decay orders, while {[}14{]} includes
the Laguerre step-function class. The present organization makes the
sharp numerator \(V(f)M_n\), its geometric interpretation, and its use
in admissible corridors explicit. This supplies a common mathematical
description for coefficient estimates and the controlled-distortion
application without asserting priority for the general inequalities.

\hypertarget{acknowledgments}{%
\subsection{Acknowledgments}\label{acknowledgments}}

This manuscript was developed in dialogue with ChatGPT (OpenAI), which
provided substantial assistance in deriving and checking family
estimates, analyzing normalization, searching for historical references,
and drafting and organizing the text. This contribution extended beyond
language editing. The investigation originated in the author's earlier
derivations of the trigonometric and Chebyshev I estimates and their use
in controlled distortion {[}1{]}; this statement describes the
development of the work, rather than a claim of priority over the
classical coefficient theory. The author directed the investigation and
examined the arguments during its development. Responsibility for
verification of the final mathematical statements and references remains
with the author.

\clearpage
\hypertarget{references}{%
\subsection{References}\label{references}}

{[}1{]} V. Sluchak. \emph{Estimations of Fourier Coefficients For
Controlled Distortion of a Periodic Function}. arXiv:2412.10376, 2024.
\url{https://doi.org/10.48550/arXiv.2412.10376}

{[}2{]} NIST Digital Library of Mathematical Functions. Chapter 18,
especially Sections 18.2, 18.3, 18.5, 18.7, 18.9, 18.14, and 18.15.
\url{https://dlmf.nist.gov/18} Accessed 17 October 2026.

{[}3{]} G. Szegő. \emph{Orthogonal Polynomials}. Fourth edition.
American Mathematical Society Colloquium Publications, Volume 23.
American Mathematical Society, Providence, RI, 1975. Classical
background for the identities and asymptotics collected in {[}2{]}.

{[}4{]} M. Taibleson. Fourier coefficients of functions of bounded
variation. \emph{Proceedings of the American Mathematical Society} 18
(1967), 766.

{[}5{]} L. Gogoladze and V. Tsagareishvili. Some classes of functions
and Fourier coefficients with respect to general orthonormal systems.
\emph{Proceedings of the Steklov Institute of Mathematics} 280 (2013),
156--168. \url{https://doi.org/10.1134/S0081543813010100}

{[}6{]} R. E. Edwards. \emph{Fourier Series: A Modern Introduction}.
Volume 1, second edition. Graduate Texts in Mathematics 64. Springer,
New York, 1979. \url{https://doi.org/10.1007/978-1-4612-6208-4}

{[}7{]} V. Sh. Tsagareishvili. On the Fourier coefficients of a function
of bounded variation {[}in Russian: O koeffitsientakh Fur'e funktsii s
konechnym izmeneniem{]}. \emph{Soobshcheniya Akademii Nauk Gruzinskoi
SSR} 101(3) (1981), 537--539.

{[}8{]} V. Tsagareishvili. On the Fourier coefficients for general
orthonormal systems. \emph{Proceedings of A. Razmadze Mathematical
Institute} 124 (2000), 131--150.
\url{https://rmi.tsu.ge/proceedings/volumes/pdf/v124-8.pdf}

{[}9{]} L. Gogoladze and V. Tsagareishvili. The Fourier coefficients of
functions with bounded variation. \emph{Russian Mathematics} 57(8)
(2013), 10--19; Russian original in \emph{Izvestiya Vysshikh Uchebnykh
Zavedenii. Matematika} (8) (2013), 14--23.
\url{https://doi.org/10.3103/S1066369X13080021}

{[}10{]} S. S. Dragomir. Approximating the Stieltjes integral of bounded
functions and applications for three point quadrature rules.
\emph{Bulletin of the Korean Mathematical Society} 44(3) (2007),
523--536. Open preprint:
\url{https://rgmia.org/papers/v10e/ASIBFApp3PQR.pdf}

{[}11{]} H. Wang. New error bounds for Legendre approximations of
differentiable functions. \emph{Journal of Fourier Analysis and
Applications} 29 (2023), Article 42. Preprint:
\url{https://arxiv.org/abs/2111.03833}

{[}12{]} L. Durand. Nicholson-type integrals for products of Gegenbauer
functions and related topics. In R. A. Askey (ed.), \emph{Theory and
Application of Special Functions}. Academic Press, 1975, pp.~353--374.

{[}13{]} S. Xiang and G. Liu. Optimal decay rates on the asymptotics of
orthogonal polynomial expansions for functions of limited regularities.
\emph{Numerische Mathematik} 145 (2020), 117--148.
\url{https://doi.org/10.1007/s00211-020-01113-3}

{[}14{]} Y. Zhang, S. Xiang, and D. Kong. On optimal convergence rates
of Laguerre polynomial expansions for piecewise functions. \emph{Journal
of Computational and Applied Mathematics} 425 (2023), Article 115053.
\url{https://doi.org/10.1016/j.cam.2022.115053}

[15] V. Sluchak. \emph{Method of control of velocity distribution harmonics forward a Propeller}. The Shipbuilding Industry, RUMB, ISS, 6s, 1986, Leningrad. (In Russian.)

[16] V. Sluchak. \emph{Flow Improving Devices For Marine Propulsers. Concept Design Study}. arXiv:2503.09616, 2025. \url{https://doi.org/10.48550/arXiv.2503.09616}

[17] V. Sluchak. \emph{Thin ring wing as a means of flow improvement upstream of a propeller}. arXiv:2502.05231, 2025. \url{https://doi.org/10.48550/arXiv.2502.05231}

\end{document}